\magnification 1200
\input amssym.def
\input amssym.tex
\parindent = 40 pt
\parskip = 12 pt
\font \heading = cmbx10 at 12 true pt
 at 22 true pt
 at 19 true pt
 at 7 true pt
\def \R{{\bf R}}
\def \C{{\bf C}}

\centerline{\heading Adapted Coordinates in Two Dimensions}
\centerline{\heading and a Proof of Puiseux's Theorem }
\rm
\line{}
\line{}
\centerline{\heading Michael Greenblatt}
\line{}
\centerline{July 24, 2008}
\baselineskip = 12 pt
\font \heading = cmbx10 at 14 true pt
\line{}
\line{}

\vfootnote{}{This research was supported in part by NSF grant DMS-0654073} 
\vfootnote{}{2000 MSC: 42B99 (primary), 32S05 (secondary)} \noindent {\bf 1. Introduction} 

Let $\C\{x_1,...,x_n\}$
denote the ring of power series whose coefficients increase slowly enough so that the series converges 
in a neighborhood of the origin in $\C^n$. Suppose $f(x,y) \in \C\{x,y\}$ with $f(0,0) = 0$. Then one 
version of Puiseux's theorem is the statement that there exists a factorization
$$f(x,y) = u(x,y)x^c \prod_{i=1}^m(y - g_i(x)) \eqno (1.1)$$
Here for some natural number $n$, each $g_i \in \C\{x^{{1 \over n}}\}$ with $g_i(0) = 0$, and $u(x,y) \in \C\{x^{{1 \over n}}
,y^{{1 \over n}}\}$ with $u(0,0) \neq 0$. Hence the zeroes of $f(x,y)$ are parameterized by analytic 
functions of one variable. (With a little more effort one can show  $u(x,y) \in \C\{x,y\}$). One
method to prove the factorization $(1.1)$ goes 
back to Isaac Newton himself. Newton's method produces the terms of the $g_i(x)$ through an infinite
recursion; in modern treatments one then shows the resulting power series converges in a 
neighborhood of the origin. The latter is normally done by invoking a topological argument involving 
Riemann surfaces (see [BK]). Alternatively, one may carefully examine the properties of 
Newton's algorithm as one proceeds and then directly prove that the 
resulting $g_i(x)$ are in some $\C\{x^{{1 \over n}}\}$; this is done in [Ca] and [Ch] 
(Puiseux's original proof was somewhat different).

The purpose of this paper is to provide an argument based on Newton's method and some ideas 
from resolution of singularities that gives a quick proof of the factorization $(1.1)$ 
(including the convergence of the $g_i(x)$). It is then shown that similar ideas can be used to give
a short proof of the existence of smooth adapted coordinates in 
two dimensions (Theorem 1.2 below). This result was first proved in the real-analytic case by 
Varchenko [V] and then recently for the general smooth case by Ikromov-M\"uller [IM]. These proofs use 
detailed information about the zero set of $S(x,y)$.

The arguments of this paper will use only the
two-dimensional implicit function theorem and some basic properties of Newton polygons; they 
are however inspired by more 
modern resolution of singularities ideas as will be discussed at the end of section 2. It should also
be pointed out that if one is willing to assume the Weierstrass preparation theorem and Hensel's lemma,
there exist short and rather different elementary proofs of Puiseux's theorem of a more algebraic 
nature. We refer to [N] for more information.

\noindent We make extensive use of the following object, essentially used in Newton's letter. 

\noindent {\bf Definition 1.2.} Let $\sum_{a,b} s_{ab}x^ay^b$ be a power series in
$x^{{1 \over n}}$ and $y^{{1 \over n}}$ for some positive integer $n$, and assume that at least 
one $s_{ab}$ is nonzero. For any $(a,b)$ for which $s_{ab} \neq 0$, let $Q_{ab}$ be the quadrant 
$\{(x,y) \in \R^2: 
x \geq a, y \geq b \}$. Then the {\it Newton polygon} $N(S)$ is defined to be 
the convex hull of all $Q_{ab}$.  

The boundary of a Newton polygon consists of finitely many (possibly zero) bounded edges 
of negative slope as well as an unbounded vertical ray and an unbounded horizontal ray. We also will
make use the following.

\noindent {\bf Definition 1.3.} Let $S(x,y)$ be as in Definition 1.2.
The {\it Newton distance} $d(S)$ is defined by $d(S) = \inf \{d: (d,d) \in N(S)\}$. 

\noindent {\bf Definition 1.4.} Suppose $e$ is a compact edge of $N(S)$. Define $S_e(x,y)$
by $S_e(x,y) = \sum_{(a,b) \in e} s_{a,b} x^ay^b$. In other words $S_e(x,y)$ is the sum of the terms
of the Taylor expansion of $S$ over all $(a,b) \in e$.

In the study of oscillatory integrals in two dimensions, the notion of adapted coordinates plays an
important role.

\noindent {\bf Definition 1.5.} A coordinate system is said to be {\it adapted} if $d(S) = \sup_{\alpha}
d(S\circ\alpha)$, where the supremum is taken over all smooth coordinate changes $\alpha$ defined in a 
neighborhood of $(0,0)$ such that $\alpha(0,0) = (0,0)$.

The significance of adapted coordinate systems is the following. Consider the oscillatory integral
$$J_{\lambda} = \int_{\R^2} e^{i \lambda S(x,y)} \phi(x,y)\,dx\,dy\eqno (1.2)$$
Assume $S(0,0) = 0$ and that $S$ has a critical point at the origin; that is, $\nabla S(0,0) = 0$. 
The function $\phi(x,y)$ is a cutoff function in a neighborhood of the origin and $\lambda$ denotes a
real parameter which one assumes to be large. Then the best (supremal) 
$\epsilon$ for which one has the estimate $|J_{\lambda}| < C\ln(2 + |\lambda|)|\lambda|^{-\epsilon}$ 
for all $\lambda$ and all $\phi(x,y)$ supported in a sufficiently small neighborhood of the origin has the nice form $\epsilon = {1 \over d(S)}$ if and only
if $S(x,y)$ is in an adapted coordinate system. This
was proved by Ikromov, Kempe, and M\"uller in [IKM]. For the real-analytic case, where one considers
real-analytic phase and takes the supremum of $(1.2)$ over real-analytic coordinate changes, the
corresponding result was earlier proven by Varchenko in [V]. 

\noindent {\bf Theorem 1.1.} Suppose $S$ has nonvanishing Taylor expansion at $(0,0)$ and $S(0,0) = 0$.
A coordinate system is adapted if any of the following three cases hold.

\noindent {\bf Case 1.} The line $y = x$ intersects $N(S)$ in the interior of a bounded edge 
$e$ and any real zero $r$ of $S_e(1,y)$ or $S_e(-1,y)$ with $r \neq 0$ has order less than $d(S)$.

\noindent {\bf Case 2.} The line $y = x$ intersects $N(S)$ at a vertex $(d,d)$.

\noindent {\bf Case 3.} The line $y = x$ intersects $N(S)$ in the interior of one of the unbounded
edges. 

\noindent {\bf Proof:} By the main theorem of [G], if $U$ is a small enough neighborhood of the origin, and $\epsilon_0$ denotes the supremum
of the numbers $\epsilon$ for which $\int_U |f|^{-\epsilon}$ is finite, then $d(S) \leq {1 \over 
\epsilon_0}$, with $d(S) = {1 \over \epsilon_0}$ in cases 1, 2, and 3. Hence if one is in cases 1, 2, 
or 3, one is in adapted coordinates.

\noindent {\bf Theorem 1.2.} Suppose $S$ has nonvanishing Taylor expansion at $(0,0)$ and 
$S(0,0) = 0$. Then there exists some coordinate system in a neighborhood of the origin such that
one of Cases 1, 2, or 3 hold. Hence there exists an adapted coordinate system for $S(x,y)$. The 
associated coordinate change can always be taken to be of the form $(x,y) \rightarrow (x, y - 
\psi(x))$ or $(x,y) \rightarrow (x - \psi(y),y)$ for a smooth $\psi$. 

In [IM], a slightly weaker version of Theorem 1.2 is proven which also shows that for any smooth phase 
there exists a 
smooth adapted coordinate system. The arguments of [IM] use Puiseux's theorem and do a careful 
analysis of the different $g_i(x)$. The proof of Theorem 1.2 is in section 3. It should be also 
pointed out that Theorem 1.1 follows from Theorem 3.3 of [IM].

\noindent {\bf 2. Proof of Puiseux's Theorem.}

Suppose $f(x,y) \in \C\{x,y\}$. After factoring out the largest possible power 
of $x$ out of $f(x,y)$, we can write $f(x,y) = x^cg(x,y)$, where $\partial_y^e g(0,0) \neq 0$ for
some $e$. Since $(1.1)$ trivially holds if $e$ is zero, we can assume $e > 0$. Assuming $e$ to be 
chosen minimal, we have that $(0,e)$ is on the Newton
polygon $N(g)$. We will prove Puiseux's theorem by proving the following theorem:

\noindent {\bf Theorem 2.1.} Suppose that $h(x,y) = \sum_{a,b} h_{ab}x^ay^b \in \C\{x^{{1 \over n}},
y\}$ such that $h(0,0) = 0$ and such that $(0,E) \in N(h)$ for some $E > 0$. Then one has a 
factorization $h(x,y) = H(x,y)(y - g(x))$ where
for some natural number $N$, $H(x,y) \in \C\{x^{{1 \over N}},y\}$ and $g \in \C\{x^{{1 \over N}}\}$
with $g(0) = 0$.

Puiseux's theorem follows by applying Theorem 2.1 repeatedly; $(0,E-1) \in N(H)$ and thus starting
with $g(x,y)$, after $e$ iterations one has $(1.1)$.

\noindent {\bf Proof of Theorem 2.1}. If $y$ divides $h(x,y)$ we are done, so we may assume that there
is some point $(D,0)$ on the Newton polygon $N(h)$ with $D > 0$. Let $(p,q)$ denote the vertex of $N(h)$ with 
$q > 0$ such that $q$ is minimal.
Thus the segment $e$ connecting $(p,q)$ to $(D,0)$ is an edge of $N(h)$. Let $h_e(x,y)$ denote
the sum of the terms $h_{ab}x^ay^b$ of the series $h(x,y) = \sum_{a,b} h_{ab}x^ay^b$ such that $(a,b)$
is on $e$. Thus $h_e(x,y)$ is a polynomial in $x^{{1 \over n}}$ and $y$. Write the equation of
the edge $e$ as $x + my = \alpha$. Hence if $h_{ab}x^ay^b$ appears in $h_e(x,y)$ then $a + mb = \alpha$.
We factor out $x^{\alpha}$, writing $h_e(x,y) = x^{\alpha}h_e'(x,y)$. Each term of $h_e'(x,y)$ is
now of the form $h_{ab}x^{a - \alpha}y^b$ and $(a - \alpha) + mb = 0$ or $(a - \alpha) = -mb$.
Thus we have
$$h_{ab}x^{a - \alpha}y^b = h_{ab}({y \over x^m})^b \eqno (2.1)$$
Conequently for a polynomial $P(z)$, we can write
$$h_e(x,y) = x^{\alpha}P({y \over x^m}) \eqno (2.2)$$
The proof of Theorem 2.1 will now proceed by an inductive process. At each stage we will perform a 
coordinate change of
the form $(x,y) \rightarrow (x,y + a(x))$ for some $a(x) \in \C\{x^{1 \over N}\}$
with $a(0) = 0$. The resulting function $h(x,y + a(x))$ will fall into one of the following two 
(not mutually exclusive) cases.

\noindent {\bf Case 1}: $y$ divides $h(x,y + a(x))$.

\noindent {\bf Case 2}: $h(x,y + a(x))$ satisfies the hypotheses of Theorem 2.1 and
the second-lowest vertex $(p'',q'')$ of the Newton polygon of $h(x,y + a(x))$ satisfies $q'' < q$. 

In the first case, one transfers back to the original coordinates and we have the conclusions of Theorem 2.1.
In the second case, one is back under the assumptions of Theorem 2.1 and thus can repeat the upcoming
argument, finding the next $a(x)$. Since $q'' < q$, after at most $q$ iterations one will have to be 
in the first case and we will be done.

So our task is to show that under the assumptions of Theorem 2.1 we can always find an $a(x)$ such that 
one of the two cases holds. Suppose first that the polynomial $P$ has a (complex) root $r$ of order $q' < q$. 
Then the function $P(z + r)$ has a root at $z = 0$ of order $q'$. Hence there is a term of 
$h_e(x,y + rx^m) = x^{\alpha}P({y \over x^m} + r)$ with $y$ appearing to the $q'$th power, but no terms
with $y$ appearing to a lower power than $q'$. Define $H(x,y) = h(x,y + rx^m)$. Note that $H_e(x,y)
= h_e(x,y + rx^m)$. Thus a segment of the line $x + my = \alpha$ is an edge of the Newton
polygon $N(H)$ of $H$, as was the case for $h$. However, instead of going down to $(D,0)$, for $H$ the
segment terminates at $(p',q')$ for some $p'$. Hence either $(p',q')$ is the lowest vertex of $N(H)$,
in which case one is in Case 1 with $a(x) = rx^m$, or the second-lowest vertex $(p'',q'')$ of $N(H)$ 
(which could be $(p',q')$) satisfies $q'' \leq q' < q$. Therefore if we let $a(x) = rx^m$, $h(x,y)$ is 
in case 2 and we are done.

Thus it remains to analyze the situation where $P(z)$ has a single complex root $r$ of order $q$. Here $P(z) = 
c(z - r)^q$ for some $c$. This is the situation where Newton's method gives an infinite iteration; here
we will do something different. We look at the function $h(x,x^m y)$. Since $x + my = \alpha$ is a 
supporting line for $N(h)$, the terms of $h_e(x,x^my)$ are the terms of $h(x,x^my)$ with with the
lowest power of $x$ appearing. Since $h_e(x,x^my) = x^{\alpha}P({x^m y \over x^m}) = $
$cx^{\alpha}(y - r)^q$, for some $\epsilon > 0$ we may write
$$h(x,x^my) = cx^{\alpha}(y - r)^q + x^{\alpha + \epsilon}l(x,y) \eqno (2.3)$$
By $(2.3)$, the function $h'(x,y) = {h(x,x^my) \over x^{\alpha}}$ is in $\C\{x^{{1 \over N}},y\}$ for
some $N$ and we have
$$h'(x,y) = c(y - r)^q + x^{\epsilon}l(x,y) \eqno (2.4)$$
The trick is now as follows. The function ${\partial^{q-1} h' \over \partial y^{q-1}}$ has a zero at
$(0,r)$, but has non-vanishing $y$ derivative there. Hence by applying the 2-dimensional implicit function
theorem (technically to  ${\partial^{q-1} h' \over \partial y^{q-1}}(x^N,y)$), one has that there
is some function $k(x) \in \C\{x^{{1 \over N}}\}$ with $k(0) = r$ such that
${\partial^{q-1} h' \over \partial y^{q-1}}(x,k(x)) = 0$ near the origin. One now defines 
$H(x,y) = h(x,y + x^mk(x))$. The fact that allows the induction to proceed is that
$${\partial^{q-1} H \over \partial y^{q-1}}(x,0) = {\partial^{q-1} h \over \partial y^{q-1}}
(x,x^mk(x)) = x^{\alpha - (q - 1)m}{\partial^{q-1} h' \over \partial y^{q-1}}(x,k(x)) = 0 \eqno 
(2.5)$$
Like before, the coordinate change is such that $x + my = \alpha$ is still a supporting line for 
$N(H)$. This time it intersects $N(H)$ in the single vertex $(p,q)$. This may be easiest to see
from $(2.3)$ using the fact that in the coordinates of $(2.3)$ the coordinate change is of the form 
$(x,y) \rightarrow (x,y + r + \tilde{k}(x))$ where $\tilde{k}(0) = 0$. 

If $y$ divides $H$ we are back
in case 1 and we are done. So we may assume there is some vertex $(d',0)$ on $N(H)$ with $d' > 0$.
If $(p,q)$ is
anything other than the second-lowest vertex, we are in case 2 and thus we'd be done again. Hence
we can assume that the segment $e'$ connecting $(p,q)$ to $(d',0)$ is an edge of $N(H)$. The condition
$(2.5)$ ensures that $H_{e'}(x,y)$ cannot have a single complex root of order $q$; for if this were to happen 
like above $H_{e'}(x,y)$ would be of the form $cx^{\alpha'}({y \over x^{m'}} - r')^q$. But this expression
has a nonvanishing $y^{q - 1}$ term; this contradicts $(2.5)$ which implies that for every $a$ the Taylor 
series coefficient $H_{a\,q-1}$ is zero. Hence $H_{e'}(x,y)$ must have a root of order less than $q$. We dealt
with this situation above; a further coordinate change of the correct form puts us in case 1 or 2. 
This completes the proof of Theorem 2.1.

Those familiar with resolution of singularities algorithms can recognize this idea of taking the
zero set of the $(q-1)$st derivative of a function and making it a hyperplane, so that an inductive 
procedure may continue. So essentially what is happening here is that an argument of this type is 
being incorporated into Newton's method to construct a process that terminates after finitely many
applications of the implicit function theorem rather than an infinite iteration.

\noindent {\bf 3. Proof of Theorem 1.2.} 

We now assume that $S(x,y)$ is a smooth function defined in a 
neighborhood of the origin with $S(0,0) = 0$ and having a nonvanishing Taylor expansion at $(0,0)$.
Let $N(S)$ denote the Newton polygon of this Taylor expansion. Assume we are not in any of the three
cases of Theorem 1.2. Thus
the line $y = x$ intersects the Newton polygon $N(S)$ in the interior of a compact edge $e$, and
$S_e(1,y)$ or $S_e(-1,y)$ has a real zero $r \neq 0$ of order $k \geq d(S)$. Replacing $x$ by $-x$ 
and/or $y$ by $-y$ if necessary, we may assume $r > 0$ is a zero of $S_e(1,y)$ of order $k$. 
The goal is to do a coordinate change of the proper form that puts
us into one of these three cases. Denote the equation of the line containing $e$ by $x + my = \alpha$. 
Exactly as $(2.2)$, there is some polynomial $Q(y)$ such that for $x > 0$ we have
$$S_e(x,y) = x^{\alpha}Q({y \over x^m})$$
Plugging in $x = 1$, we see that $Q(y) = S_e(1,y)$. Hence $S_e(x,y)$ has zeroes 
of order $k$  on the curve $y = r x^m$. This implies that $S_e(x,1)$ has a zero of order $k$ at 
$x = r^{-{1 \over m}}$. As a result, we may switch the roles
of the $x$ and $y$ axes if necessary and assume that $m \geq 1$; this makes our subsequent arguments
somewhat easier. 

Next, we show that $m$ must in fact be an integer. To see this, note that if $m$ were not an integer, 
then the degrees of the powers of $y$
appearing in $S_e(1,y)$ would have to be separated by at least 2. Hence $S_e(1,y)$ would have to be 
of the form $y^{\beta}R(y^c)$ for some $\beta \geq 0,$ $c \geq 2$, where $R$ is a polynomial.  Next, since 
$(d(S),d(S))$ is
on $N(S)$, we have $\alpha = (1 + m)d(S)$. Since  $m > 1$ when $m \geq 1$ is not an integer, the 
maximum possible value of $y$ on the line $x + my = (1 + m)d(S)$ for $x, y \geq 0$ is 
${m + 1 \over m} d(S) < 2d(S)$. Thus the degree of $y^{\beta}R(y^c)$ is less than $2d(S)$, and
hence the degree of $R(y)$ is less than ${2d(S) \over c} \leq d(S)$. Hence the zeroes of $R(y)$ are
of order less than $d(S)$, implying the zeroes of 
of $S_e(1,y) = y^{\beta}R(y^c)$ other than $y = 0$ are of order less than $d(S)$. This
contradicts our assumption that $S_e(1,y)$ has a zero of order $k \geq d(S)$ and we conclude that $m$ is an
integer. 

Note that if $m$ is even, then $S_e(1,y) = \pm S_e(-1,y)$, while if $m$ is odd one
has $S_e(1,y) = \pm S_e(-1,-y)$. Hence both $S_e(1,y)$ and $S_e(-1,y)$ have a zero of order $k$; we 
never really had to replace $x$ by $-x$ in the above. The preceding arguments thus imply:

\noindent {\bf Fact:} If one is not in adapted coordinates and $m \geq 1$, then $m$ is an integer
and $Q(y) = S_e(1,y)$ has a zero of order $k \geq d(S)$.

We now come to the main argument; we will prove the existence of a coordinate change of the 
form $(x,y) \rightarrow (x,y + a(x))$, $a(x)$ smooth, that puts us into one of the three cases. (The coordinate change
$(x,y) \rightarrow (x + a(y),y)$ corresponds to $m \leq 1$). We proceed as follows. Let $(p,q)$ denote
the upper vertex of
the edge $e$; necessarily $q > d(S)$. We will find a smooth function $a(x)$ such that  
$S'(x,y) = S(x,y + a(x))$ is in one of the following two (not mutually exclusive) categories.

\noindent {\bf Category 1:} $S'(x,y)$ is in one of the three cases of adapted coordinates.

\noindent {\bf Category 2:} The line $y = x$ intersects the interior of an edge $e'$ of $N(S')$
with equation $x + m'y = \alpha '$, $m' \geq 1$, such that the upper vertex $(p',q')$ of $e'$ 
satisfies $q' < q$.

\noindent Theorem 1.2 will then follow; there can be at most $q$ iterations of category 2.

The arguments now resemble those of section 2. We first consider the case where the order $k$ of the
zero $r$ of 
$Q(y) = S_e(1,y)$ satisfies $k < q$. Here we let $a(x) = rx^m$, and thus $S'(x,y) = S(x,y + rx^m)$. Then 
$x + my = \alpha$ is a supporting line of $N(S')$ as it was for $N(S)$, and like in section $2$ there is
an edge $E$ of $N(S')$ on this line whose upper vertex is $(p,q)$. Note that $S'_E(x,y) = S_e(x,y + 
rx^m) = x^{\alpha}Q({y \over x^m} + r)$. Since $Q$ has a zero of order $k$ at $r$, the lowest power of 
$y$ appearing in $S'_E(x,y)$ is $y^k$, and 
therefore $E$'s lower vertex is at a point $(j,k)$ for some $j$. Since both vertices of $E$ have 
$y$-coordinates
at least $d(S)$, they are both in the portion of $x + my = \alpha$ on or above $(d(S),d(S))$. Thus
the edge $E$ lies wholly on or above the line $y = x$. If the line $y = x$ intersects
$N(S')$ at a vertex or inside the horizontal or vertical rays, one is in Category 1. Otherwise, 
it must intersect $N(S')$ in the 
interior of an edge $e'$ whose upper vertex is either $(j,k)$ or a lower vertex. And because
$x + my = \alpha$ is a supporting line for $N(S')$ and $e'$ lies below $E$, $e'$ will have equation 
$x + m'y = \alpha '$ for some $m' > m \geq 1$. Thus we are in category 2. Hence when $k < q$ we are in 
either Category 1 or 2 and we are done.

It remains to consider the situation where $r$ is a zero of $Q(y)$ of order $q$. In this case we have
$Q(y) = c(y - r)^q$ for some $c$. For a large integer $n$ we expand $S(x,y)$ as
$$S(x,y) = cx^{\alpha}({y \over x^m} - r)^q + T_n(x,y) + E_n(x,y) \eqno (3.1)$$
Here the polynomial $T_n(x,y)$ are the terms of $S$'s Taylor expansion with exponents less
than $n$. For all $0 \leq \beta, \gamma < n$ one has
$$ |{\partial^{\beta + \gamma}E_n \over \partial x^{\beta} \partial y^{\gamma}}(x,y)| < C(|x|^{n - \beta} + 
|y|^{n - \gamma}) \eqno (3.2)$$
Analogous to $(2.3)$, one has
$$S(x,x^my) =cx^{\alpha}(y - r)^q + x^{\alpha + 1}T_n'(x,y) + E_n(x,x^my) \eqno
(3.3)$$
Here $T_n'(x,y)$ is also a polynomial. Analogous to $(2.4)$ we define
$s(x,y) = {S(x,x^my) \over x^{\alpha}}$, so that
$$s(x,y) = c(y - r)^q + xT_n'(x,y) + x^{-\alpha}E_n(x,x^my) \eqno (3.4)$$
We claim that the function $s(x,y)$ is smooth on a neighborhood of $(0,r)$. Off the $y$-axis smoothness
holds because $S(x,y)$ is smooth. One can show that a given derivative of $s(x,y)$ exists when
$x = 0$ and equals that of $c(y - r)^q + xT_n'(x,y)$ for large enough $n$ by
examining the difference quotient of a one-lower order derivative of $(3.4)$, inductively assuming 
this lower-order derivative exists and has the right value when $x = 0$. Equation
$(3.2)$ ensures that the difference quotient of the lower derivative of $x^{-\alpha}E_n(x,x^my)$ tends 
to zero as $x$ goes to zero.
We conclude that $s(x,y)$ is smooth on a neighborhood of $(0,r)$. 

Analogous to after $(2.4)$, we next use the smooth implicit function 
theorem on ${\partial^{q-1}s \over \partial y^{q-1}}$ and find a smooth function $k(x)$ defined
in a neighborhood of $x = 0$ such that $k(0) = r$ and
${\partial^{q-1}s \over \partial y^{q-1}} (x,k(x)) = 0$.
Transferring this back to $S(x,y)$, as in $(2.5)$ we have
$${\partial^{q-1}S \over \partial y^{q-1}} (x,x^mk(x)) = 0 \eqno (3.5)$$
Thus if we let $a(x) = x^mk(x)$ and $S'(x,y) = S(x,y + x^mk(x))$, for all $x$ we consequently have
$${\partial^{q-1}S' \over \partial y^{q-1}} (x,0) = 0 \eqno (3.6)$$
Thus for every $a$ the Taylor series coefficient $S'_{a\,q-1}$ is zero.

Next, since $x + my = 
\alpha$ is a supporting line for $N(S)$, analogous to after $(2.5)$ this line is also a supporting
line for $N(S')$ and intersects $N(S')$ at the single vertex $(p,q)$. If $S'(x,y)$ is in adapted
coordinates we are in Category 1 and have nothing to prove, so we may assume the coordinates are not 
adapted. Let $e'$ denote the
edge of $N(S')$ intersecting the line $y = x$ and denote its equation by $x + m'y = \alpha '$. If the 
upper vertex $(p',q')$ of $N(S')$ satisfies 
$q' < q$, one is in Category 2 and we are done. So we assume this upper vertex is $(p,q)$ itself.
Also, since $e'$ lies within the set $x + my \geq \alpha$
and is no higher than the vertex $(p,q)$ of $N(S')$ that is on the supporting line $x + my = \alpha$, 
we have $m' > m \geq 1$. 

If $S'_{e'}(1,y)$ has a real zero $r' \neq 0$ of order 
less than $q$, one is in the situation above $(3.1)$; there is a smooth $b(x)$ such that 
$S'(x,y + b(x)) = S(x,y + a(x) + b(x))$ is in Category 1 or 2 as needed. The only other 
possibility is that $S'_{e'}(1,y)$ has a single zero $r' \neq 0$ of order $q$. But like at the end of section 2,
this cannot happen. For this would imply $S'_{e'}(x,y) = c'x^{\alpha '}({y \over x^{m'}}- r')^q$ 
has a nonvanishing $y^{q-1}$ term. Consequently, for some $a$ the Taylor series
coefficient $S'_{a\,q-1}$ would be nonzero, contradicting $(3.6)$. Thus the case where 
$S'_{e'}(1,y)$ has a single zero of order $q$ does not occur, and we are done.

\noindent {\bf 4. References}

\noindent [BK] E. Brieskorn and H. Knorrer, {\it Plane Algebraic Curves}, Birkhauser Verlag, Basel, 
(1986), vi+721 pp. \parskip = 3pt

\noindent [Ca] E. Casas-Alvero, {\it Singularities of plane curves}, London Mathematical Society 
Lecture Note Series {\bf 276} (2000), Cambridge University Press, Cambridge, xvi+345 pp.

\noindent [Ch] A. Chenciner, {\it Courbes alg\'ebriques planes}, Publications Math\'ematiques de 
l'Universit\'e Paris VII {\bf 4} (1978), Universit\'e de Paris VII, U.E.R. de Math\'ematiques, Paris,
203 pp.

\noindent [G] M. Greenblatt, {\it Newton polygons and local integrability of negative powers of smooth 
functions in the plane}, Trans. Amer. Math. Soc. {\bf 358} (2006), no. 2, 657-670.

\noindent [IM] I. Ikromov and D. M\"uller, {\it On adapted coordinate systems}, preprint.

\noindent [IKM] I. Ikromov, M. Kempe, and D. M\"uller, {\it Sharp $L^p$ estimates for maximal operators
associated to hypersurfaces in $\R^3$ for $p > 2$}, preprint.

\noindent [N] K. Nowak, {\it Some elementary proofs of Puiseux's theorems}, Univ. Iagel. Acta Math. 
{\bf 38} (2000), 279-282.

\noindent [V] A. N. Varchenko, {\it Newton polyhedra and estimates of oscillatory integrals}, 
Functional Anal. Appl. {\bf 18} (1976), no. 3, 175-196.

\line{}
\line{}

\noindent Department of Mathematics, Statistics, and Computer Science \hfill \break
\noindent University of Illinois at Chicago \hfill \break
\noindent 322 Science and Engineering Offices \hfill \break
\noindent 851 S. Morgan Street \hfill \break
\noindent Chicago, IL 60607-7045 \hfill \break

\end